\documentclass[12pt,a4paper]{article}
\usepackage{amsmath, amssymb, theorem, latexsym, epsfig}
 \usepackage{graphics}
\usepackage{epstopdf}

\newtheorem{theorem}{Theorem}[section]
\newtheorem{lemma}[theorem]{Lemma}
\newtheorem{proposition}[theorem]{Proposition}
\newtheorem{definition}[theorem]{Definition}

\newtheorem{conjecture}[theorem]{Conjecture}
\newtheorem{corollary}[theorem]{Corollary}
\newtheorem{example}[theorem]{Example}

\allowdisplaybreaks

\def\neweq#1{\begin{equation}\label{#1}}
\def\endeq{\end{equation}}
\numberwithin{equation}{section}
\begin{document}
{\centering
\bfseries
{\LARGE Even Subdivison-Factors of\\ 
\vspace{0.3cm}
 Cubic Graphs }\\
 
\bigskip
\mdseries
{\large Arthur Hoffmann-Ostenhof}\footnote{supported by the FWF project P20543.}
\par
\upshape
\footnotesize \textit{Technical University of Vienna, Austria}\\

}

\noindent
\begin{abstract}
\noindent
We call a set $\mathcal S$ of graphs an "even subdivison-factor" of a cubic graph $G$ if $G$ contains a spanning subgraph $H$ 
such that every component of $H$ has an even number of vertices and is a subdivision of an element of $\mathcal S$. 
We show that any set of $2$-connected graphs which is an even subdivison-factor of every $3$-connected cubic graph, 
satisfies certain properties. As a consequence, we disprove a conjecture which was stated in an attempt to solve the circuit double cover conjecture. 
\end{abstract}

\noindent
Keywords: circuit double cover, factor, frame, Petersen graph 

\vspace{0.5cm}

\noindent


\section{Basic definitions and main results}
 
For terminology not defined here we refere to \cite{Bo}. 
There are several ways to describe that a spanning subgraph with certain properties exists in a cubic graph $G$.  \\ 
A set $S$ of graphs is called a \textit{component-factor} of $G$ if $G$ has a spanning subgraph $H$ such that every component of $H$ is an element of $S$, see \cite{Plum}. Within the topic of circuit double covers the notion of a \textit{frame} was introduced, see \cite{G,HM,Z1,Z2}. Some slightly different definitions of a frame exist. Here, a frame of $G$ is a graph $F$ where every component of $F$ is either an even circuit or a $2$-connected cubic graph such that the following holds: $G$ has a spanning subgraph $F'$ which is a subdivision of $F$ and every component of $F'$ has an even number of vertices. For our purpose it is useful to join these two concepts. 

\begin{definition}
A set $\mathcal S$ of graphs is called a \textit{subdivison-factor} of a cubic graph $G$ 
if $G$ contains a spanning subgraph $H$ such that every component of $H$ is a subdivision of an element of $\mathcal S$. If every component of $H$ has an even number of vertices then $\mathcal S$ is called an even subdivision-factor of $G$.
\end{definition}



\begin{example}
Every $3$-edge colorable cubic graph $G_3$ has a spanning subgraph consisting of even circuits, i.e. an even $2$-factor. Hence, $\{C_2\}$ where $C_2$ denotes the circuit of length $2$, is an even subdivision-factor of $G_3$. Reversely, if $\{C_2\}$ is an even subdivision-factor of a cubic graph $G$, it follows that $G$ is $3$-edge colorable.
\end{example}

\noindent
Thus an even subdivision-factor is a generalization of an even $2$-factor.
It was asked in a preprint of \cite {HM} whether $\{C_2\} \cup \mathcal H$ where $\mathcal H$ is a certain infinite family of hamiltonian cubic graphs, is an even subdivision-factor of every $3$-connected cubic graph. In particular the following is conjectured in \cite{HM}. 
 (A cubic graph $G$ which admits a $3$-edge coloring such that each pair of color classes forms an hamiltonian circuit, is called a \textit{Kotzig graph}, see \cite{HM,Z1}.)


\begin{conjecture}\label{3cf}
Every 3-connected cubic graph has a spanning subgraph which is a subdivision of a Kotzig graph. 
\end{conjecture}

\noindent
A positive answer to this conjecture would have solved the circuit double cover conjecture (CDCC), see \cite{HM}. 
 For stating the main theorem which provides a negative answer to Conjecture \ref{3cf} and the posed question above, we use two definitions. 
 

\begin{definition}\label{98}
Let $H_i$, $i \in \{1,2\}$ be a subgraph of a graph $G$ or a subset of $V(G)$. Denote by $[H_1,H_2]$ the set of all paths with connect a vertex of $H_1$ with a vertex of $H_2$. Then, $d_G(H_1,H_2)$ or in short $d(H_1,H_2):= \min\limits_{\alpha \in [H_1,H_2]} \,|E(\alpha)|$. %
\end{definition}

\noindent
The parameter $l(G)$ below measures to which extend $G$ is not hamiltonian.

\begin{definition}\label{}
Let $G$ be a $2$-connected graph. Denote by $U(G)$ the set of all circuits of $G$. Define 
$$l(G):= \min\limits_{C \in U(G)} \,\max\limits_{v \in V(G)} d(C,v)\,.$$ 
Let $\mathcal S$ be a set of 2-connected graphs. Define $l_m(\mathcal S):= \max\limits_{G \in \mathcal S} \,l(G)$ if this maximum exists; otherwise set 
$l_m(\mathcal S):= \infty$.
\end{definition}

\noindent
Note that in the case of $G$ being hamiltonian, $l(G)=0$.
We state the main result.

\begin{theorem}\label{infinit}
Let $\mathcal S$ be a set of $2$-connected graphs which is an even subdivision-factor of every $3$-connected cubic graph, then $l_m(\mathcal S)= \infty $.
\end{theorem}

\noindent
Theorem \ref{infinit} implies that there is no finite set of graphs which is an even subdivision factor of every $3$-connected cubic graph. 
Note that Conjecture \ref{3cf} remains open for cyclically $4$-edge connected cubic graphs. A positive answer to this version would still solve the CDCC since a minimal counterexample to the CDCC is at least cyclically $4$-edge connected. In order to prove Theorem \ref{infinit}, we prove Theorem \ref{pktk} which concerns the iterated Petersen graph. From now on, we make preparations for the proof of Theorem \ref{pktk}.

\section{The iterated Petersen graph}

\noindent
We denote by $P_{10}$ the Petersen graph and we set $P:=P_{10}-z$, $z \in V(P_{10})$.
The \textit{iterated Petersen graph} which is defined next has already been introduced in \cite{BoS}.

\begin{definition}\label{Gk}
Let $G$ be a graph with $d(v) \in \{2,3\}$, $\forall v \in V(G)$. A $P$-inflation at $v_0 \in V(G)$ is defined as the following operation:
add $P$ to $G-v_0$ and connect each former neighbor of $v_0$ to one distinct $2$-valent vertex of $P$.\\
$G^{0}, G^{1}, G^{2},...,G^{k}$ with $k\in \mathbb{N}$ and $G^{0}:=G$, is the sequence of graphs where $G^{i}$, $i\in \{1,2,...,k\}$ results from $G^{i-1}$ by applying the $P$-inflation at every vertex in $G^{i-1}$. 
We call $P_{10}^k$ for $k \geq 1$ an iterated Petersen graph.
\end{definition}

\begin{figure}[htpb] 
\centering\epsfig{file=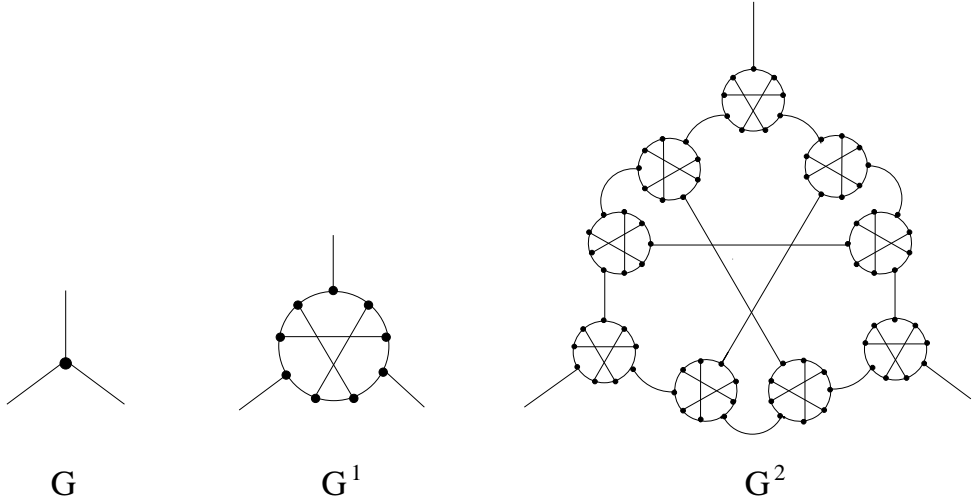,width=5in}
\caption{A vertex in a cubic graph $G$ and the corresponding copies of $P^{i-1}$ in $G^i$, $i=1,2$.}\label{tripp}
\end{figure}

\noindent
Obviously, $G^k$ is cubic if $G$ is cubic. If $G$ is not cubic, then $G$ and $G^k$ have the same number of vertices of degree $2$. 
See Figure \ref{tripp} for an illustration of Def. \ref{Gk}. Note that if we remove in the illustration of $G^i$ the dangling edges, we obtain $P^{i-1}$, $i=1,2$.

\begin{definition}\label{wkdk}
 Let $W_k$, $k \in \mathbb N$ denote the set of the three $2$-valent vertices of $P^k$ and set $d_k:= \max \,\{\, d(W_k,v) \,|\, v \in V(P^k) \,\}$.
 If a graph $X$, say, is isomorphic to $P^k$, then $W_k(X)$ denotes the set of the three $2$-valent vertices.
\end{definition}

\begin{proposition}\label{dk}
Let $k \in \mathbb N$, then $d_k= 2^{2k+1}-1$.
\end{proposition}

\noindent
Proof: The statement obviously holds for $k=0$. Consider $P^k$ for $k>0$ and set $j_k:= \min \,\{\,|V(\alpha)|\,\,|\,\,\alpha \in [w_1,w_2]\,\}$ with $\{w_1,w_2\} \subseteq W_k$ and $w_1 \not= w_2$. Let $k \geq 1$, then $P^{k}$ contains $9$ disjoint copies of $P^{k-1}$. 
$P$ results from $P^k$ by contracting each of them to a distinct vertex. Hence, every copy $P'$ of $P^{k-1}$ in $P^k$ corresponds to a vertex in $P$. 
We say a path $\alpha$ \textit{traverses} $P' \subseteq P^k$ if $\alpha$ contains a subpath $\alpha ' \subseteq P'$ which connects two distinct vertices of 
$W_{k-1}(P')$.
Every shortest path in $P^k$ which connects $w_1$ with $w_2$, traverses exactly $4$ copies of $P^{k-1}$ and thus $j_k = 4 \,j_{k-1}\,$. Since $j_0=4$, we obtain $$j_k=4^{k+1}\,.\,\,\,\,\,\,\,\,\,\,\,\,\,\,\,\,\,\,\,\,\,\,\,\,\,\,\,\,\,\,\,\,\,\,\,\,(1)$$

\noindent
Let $k \in \mathbb N$. Set $b_k:= \max\limits_{v \in V(P^k)}\, d(w_1,v)$ and $B_k:= \{v \in V(P^k)\,|\,d(v,w_1)=b_k\}$. We claim that 
$$B_k = W_k-\{w_1\}\,.\,\,\,\,\,\,\,\,\,\,\,\,\,\,\,\,\,(2)$$
We proceed by induction on $k$. For $k=0$, the statement holds. 
Let $P' \subseteq P^k$ be a copy of $P^{k-1}$ with $v_0 \in B_k \cap V(P')$.  
Then obviously $P'$ corresponds to a $2$-valent vertex of $P$. 
Let $q_1$, $q_2$ denote the two distinct vertices of $P'$ which form together a vertex cut of $P^k$ and which are both contained in $W_{k-1}(P')$. Then, $d(w_1,q_1)=d(w_1,q_2)$. The induction assumption for $k-1$ on $P'$ implies that $v_0 \in W_{k-1}(P')$. Since $v_0 \not \in \{ q_1,q_2\}$, $v_0 \in W_k-\{w_1\}$. Hence the claim is proven.

\noindent
Let $k \geq 1$ and let now $P' \subseteq P^k$ be a copy of $P^{k-1}$ with $x \in V(P')$ and $d(x,W_k)= d_k$, see Def. \ref{wkdk}. Obviously, $P'$ corresponds to a vertex of degree $3$ in $P$.
Let $\alpha_x \subseteq P^k$ connect $x$ with a vertex of $W_k$ and satisfy $|E(\alpha_x)|=d_k$. 
Hence $\alpha_x$ is a shortest path and traverses exactly one copy of $P^{k-1}$ which corresponds to a $2$-valent vertex of $P$.  
By applying (2) on $P'$ we conclude that $x \in W_{k-1}(P')$. Thus, $|E(\alpha_x)|=2j_{k-1}-1$ which finishes the proof.

\begin{corollary}\label{Pkdk}
$l(P_{10})=1$ and $l(P^k_{10})= 2^{2k-1}$, $\forall \, k \geq 1$.
\end{corollary}

\noindent
Proof: Since $P_{10}$ has no hamiltonian circuit but $P_{10}-v_0$ is hamiltonian for every $v_0 \in V(P_{10})$, $l(P_{10})=1$. 
Let $k \geq 1$, then $P^k_{10}$ contains ten disjoint copies of $P^{k-1}$ which we denote by $X_i$, $i=1,2,...,10$.
Ever circuit in $P^k_{10}$ is vertex-disjoint with at least one $X_i$ since otherwise it would imply that $P_{10}$ is hamiltonian. 
Hence, $l(P^k_{10}) \geq d_{k-1}+1$. 
It is not difficult to see that $P^k_{10}$ contains a circuit $C$ which passes through $X_i$ for $i=1,2,...,9$ and satisfies 
$W_{k-1}(X_i) \subseteq V(C)$. By the properties of $C$ and since $ \bigcup_{i=1}^{10} {V(X_i)}=V(P^k_{10}) $, it follows that 
$d(C,v) \leq d_{k-1}+1$, $\forall v \in V(P^k_{10})$. 
Hence, $l(P^k_{10}) = d_{k-1}+1$ and by applying Prop. \ref{dk}, the proof is finished.

\subsection{f-matchings and P-inflations}

\begin{definition}\label{fm}
A matching $M$ of a cubic graph $G$ is called an $f$-matching if every component of $G-M$ is $2$-connected and has an even number of vertices. 
\end{definition}

\begin{lemma}\label{3cut2}
Suppose a cubic graph $G$ has a minimal $3$-edge cut $E_0$. Then for every $f$-matching $M$ of $G$, $|M \cap E_0| \in \{0,1\}$. 
\end{lemma}

\noindent
Proof: Suppose $|M \cap E_0|=3$. Since $E_0$ is a minimal edge-cut, $G-E_0$ consists of two components which have both an odd number of vertices. Let $L$ be one of them. Then $L-M$ and thus $G-M$ contains at least one component which has an odd number of vertices, in contradiction to Def. \ref{fm}. \\
Suppose $|M \cap E_0|=2$. Then the one edge of $E_0$ which is not contained in $M$ is a bridge in $G-M$ which contradicts Def. \ref{fm}. Hence the proof is finished.

\begin{lemma}\label{connp2}
Let $E_0:= \{e_1,e_2,e_3\}$ be a minimal $3$-edge cut in a $2$-connected cubic graph $G$ such that $P$ is one component of $G-E_0$. 
Then for every $f$-matching $M$ of $G$ the following is true.
 
\noindent
$(1)$ Consider $P \subseteq G$ as a graph and $M$ restricted to $P$. Then $P-M$ is connected.

\noindent
$(2)$ $G-M$ contains a $3$-valent vertex within $V(P)$, i.e. at least one vertex of $P \subseteq G$ is not matched by $M$.
\end{lemma}

\noindent
Proof: Let $W_0:= \{w_1,w_2,w_3\}$ denote the set of the $2$-valent vertices of $P$ and let $e_i \in E_0$ be incident with $w_i$, $i=1,2,3$. By Lemma \ref{3cut2}, $|M \cap E_0| \in \{0,1\}$.

\noindent
Proof of the first statement:

\noindent
Case 1. $|M \cap E_0|=0$.\\
All $w_i$'s are contained in the same component $L$, say, of $P-M$ since otherwise one component of $G-M$ would have $e_i$, for some $i \in \{1,2,3\}$ as a bridge in contradiction to Def. \ref{fm}. Suppose by contradiction that $P-M$ has another component $L'$. 
Since $V(L') \cap W_0= \emptyset$, $L'$ is not only a component of $P-M$ but also of $G-M$. By Def. \ref{fm}, $L'\subseteq P$ is $2$-connected and thus contains a circuit. There is exactly one circuit $C'$ in $P$ which contains no vertex of $W_0$, see Figure \ref{tripp}. Then $e_i$, $i=1,2,3$ is a bridge in $G-M$ contradicting Def. \ref{fm}. Hence $P-M$ is connected. 

\noindent
Case 2. $|M \cap E_0|=1$. Let w.l.o.g. $M \cap E_0= \{e_3\}$. \\
Then $w_1$ and $w_2$ are contained in the same component $L$, say, of $P-M$ otherwise $e_i$, $i \in \{1,2\}$ is a bridge of $G-M$. 
Suppose by contradiction that $P-M$ has another component $L'$. Since $e_3$ is matched and $w_i \in V(L)$, $i=1,2$, $L'$ is not only a component of $P-M$ but also of $G-M$. By Def. \ref{fm}, $L'$ is $2$-connected and thus contains a circuit $C'$. Since $L$ is a component, $L$ contains a path $\beta$ (which is vertex-disjoint with $C'$) connecting $w_1$ with $w_2$. $P_{10}$ is obtained from $P$ and $E_0$ by identifying the three endvertices of $e_i$, $i=1,2,3$ which are not in $P$. Then $\beta$ and $C'$ correspond to two disjoint circuits in $P_{10}$ which form a $2$-factor of 
$P_{10}$. Hence $C'=L'$, and $L'$ is a circuit of length $5$ which contradicts Def. \ref{fm}.

\noindent
Proof of the second statement:

\noindent
Suppose by contradiction that every vertex of $P$ is matched by $M$. Since $|V(P)|$ is odd and by Lemma \ref{3cut2}, $|E_0 \cap M|=1$. Such matching $M$ covering $V(P)$ corresponds to a perfect matching of $P_{10}$. Hence, $P-M$ consists of a path and a circuit $C$ of length $5$. Then $C$ is also a component of $G-M$ which contradicts Def. \ref{fm}.

\begin{lemma}\label{connp6}
Let $G$, $E_0$ and $P$ be as in the previous lemma. Let $\alpha$ be a path in $G$ which passes through $P$, i.e. $\alpha$ has no endvertex in $P$ 
and $|E(\alpha) \cap E_0|=2$. Then for every $f$-matching $M$ with $E(\alpha) \cap M = \emptyset$
the following is true: $G-M$ contains a $3$-valent vertex within $V(\alpha) \cap V(P)$, i.e. at least one vertex of $V(\alpha) \cap V(P)$ is not matched by $M$.
\end{lemma}

\noindent
Proof: Suppose by contradiction that every vertex of $V(\alpha) \cap V(P)$ is matched by $M$. Then $\alpha \cap P$ is a component of $P-M$ and thus by Lemma \ref{connp2} (1) the only component of $P-M$. Since $\alpha \cap P$ contains no $3$-valent vertex we obtain a contradiction to Lemma \ref{connp2} (2) which finishes the proof.


\begin{proposition}\label{inflationcon}
Let $G$ be a $2$-connected cubic graph and $v_0 \in V(G)$. Denote by $G'$ the cubic graph which is obtained from $G$ by applying the $P$-inflation at $v_0$. 
Then $G'-M'$ is $2$-connected for every $f$-matching $M'$ of $G'$ if and only if $G-M$ is $2$-connected for every $f$-matching $M$ of $G$.
\end{proposition}

\noindent
Proof: Denote by $P'$ the subgraph of $G'$ which is isomorphic to $P$ and corresponds to $v_0 \in V(G)$.\\
Suppose by contradiction that $M'$ is an $f$-matching of $G'$ such that $G'-M'$ is not $2$-connected whereas $G-M$ is $2$-connected for every $f$-matching $M$ of $G$. Set $M'_1:= \{ e \in M' \,|\, e \not\in E(P')\}$. Denote by $M_1$ the subset of $E(G)$ which corresponds to $M'_1$. Then, 
$$(G'-M')/V(P') = G-M_1 \,\,\,\,\,\,\,\,\,\,\,\,\,\,\,\,\,\,\,\,\,\,(3)\,$$ 
\noindent
We show that $M_1$ is an $f$-matching. Lemma \ref{3cut2} implies that $v_0 \in V(G)$ is covered by at most one edge of $M_1$. Hence, $M_1$ is a matching of $G$. Since $P'-M'$ is connected by Lemma \ref{connp2} (1), equation $(3)$ implies that $G-M_1$ has the same number of components as $G'-M'$. 
Contracting an edge or shrinking a subset of vertices in a bridgeless graph does not create a bridge. Therefore and since $G'-M'$ is bridgeless by 
Def. \ref{fm}, equation $(3)$ implies that $G-M_1$ is bridgeless. Every component of $G-M_1$ has a corresponding isomorphic component in $G'-M'$ (and thus an even number of vertices) with the one exception of the component $L_0$, say, which contains $v_0$.  
 $P'-M'$ is connected by Lemma \ref{connp2} (1). Denote by $L'_0$ the component of $G'-M'$ with $(P'-M') \subseteq L'_0$. $V(L'_0)$ differs from $V(L_0)$ by containing the vertices of $V(P'-M')$ instead of $v_0$. Since $|V(L'_0)|$ is even by Def. \ref{fm} and both $|V(P'-M')|$ and $|\{v_0\}|$ are odd, $|V(L_0)|$ is even. Hence $M_1$ is an $f$-matching of $G$. Since $G-M_1$ is not $2$-connected we obtain a contradiction to the assumption in the beginning.\\
Suppose by contradiction that $M$ is an $f$-matching of $G$ such that $G-M$ is not $2$-connected whereas $G'-M'$ is $2$-connected for every $f$-matching $M'$ of $G'$. Denote by $M'_2$ the matching of $G'$ which corresponds to $M$ of $G$ with $E(P') \cap M'_2 = \emptyset$. Then $M'_2$ is an $f$-matching of $G'$. Since $G'-M'_2$ is not $2$-connected we obtain a contradiction which finishes the proof.

\begin{corollary}\label{qq}
 For every $f$-matching $M$ of $P^k_{10}$, $k \in \mathbb{N}$, $P^k_{10}-M$ is homeomorphic to a $2$-connected cubic graph.
\end{corollary}

\noindent
Proof: $P_{10}-M$ is not a circuit since it would imply that $P_{10}$ is hamiltonian. Therefore and since every bridgeless disconnected subgraph of $P_{10}$ consists of two circuits of length $5$, $P_{10}-M$ is homeomorphic to a $2$-connected cubic graph. Since $P^k_{10}$ is not hamiltonian and results from $P_{10}$ by $P$-inflations and since Proposition \ref{inflationcon} can be applied after each $P$-inflation, the corollary follows.

\subsection{Frames}

\begin{lemma}\label{pkpkpk}
Let  $k \in \mathbb N$, then $P^{k}_{10}$ is a frame of $P^{k+1}_{10}$. 
\end{lemma}

\begin{figure}[htpb]
\centering\epsfig{file=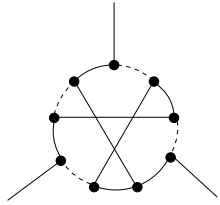,width=1.1in}
\caption{A matching of a copy of $P$ in $P^{k+1}$.} \label{Pmatch}
\end{figure}

\noindent
Proof: Let $M$ be a matching of $P^{k+1}_{10}$ such that every copy of $P$ in $P^{k+1}_{10}$ is matched as in Figure \ref{Pmatch}; $M$ is illustrated by dashed lines. 
Then $M$ is an $f$-matching of $P^{k+1}_{10}$ and the cubic graph homeomorphic to $P^{k+1}_{10}-M$ is $P^{k}_{10}$. Hence $P^{k}_{10}$ is a frame of $P^{k+1}_{10}$.

\begin{definition}\label{l2}
Let $ \alpha$ be a path in a graph $G$, then $p(\alpha)$ denotes the number of distinct copies of $P$ with which $\alpha$ has a non-empty vertex-intersection. 
For $H_i \subseteq G$, $i=1,2$, we define $p[H_1,H_2]:=\min \,\{p(\alpha) \,|\, \alpha \in [H_1,H_2]\}$
and we set $p_k:= \max \,\{p[v,W_k]\,\,|\, v \in V(P^k)\}$, $k \in \mathbb N$.
\end{definition}

\begin{lemma}\label{dkak}
Let $k \in \mathbb N$, then $p_{k+1}= 2^{2k+1}$ and $p_0=1$.
\end{lemma}

\noindent
Proof: Clearly, $p_0=1$. Let $P(x)$ and $P(y)$ denote two distinct copies of $P$ in $P^{k+1}_{10} $, $k \in \mathbb N$ with $x \in V(P(x))$ and $y \in V(P(y))$. Let $x'$ ($y'$) be the vertex in $P^k_{10}$ which corresponds to $P(x)$ ($P(y)$) by regarding $P^k_{10}$ as the graph which is obtained from $P^{k+1}_{10}$ by contracting every copy of $P$. Then for every path $\alpha \in [x,y]$ and its corresponding path $\alpha ' \in [x',y']$, $p(\alpha)=|V(\alpha ')|$. 
Hence, $p[x,y] \geq d(x',y')+1$. Since for every given path $\beta ' \in [x',y']$, there is a path $\beta \in [x,y]$ with $p(\beta)=|V(\beta')|$, $p[x,y]=d(x',y')+1$. Therefore, $p_{k+1}=d_k + 1$ (Def. \ref{wkdk}) and by applying Prop. \ref{dk} the proof is finished.


\begin{theorem}\label{pktk}
Let $\mathcal F(k)$ be the set of frames of $P_{10}^k$, $k\in \mathbb{N}$, then

\noindent
(1) every frame $G$ of $P_{10}^k$ is cubic and $2$-connected, and\\
\\
(2) $\min\limits_{G \,\in \mathcal F(k)} l(G)= \begin{cases}
  k  & \textrm {for} \,\,\,\, k \in \{ 0,1 \}\\
  2^{2k-3} & \text{for} \,\,\,\, k \geq 2 \,\,\,\,\,\,\,. 
  \end{cases}$


\end{theorem}

\noindent
Proof: Corollary \ref{qq} implies that every element of $\mathcal F(k)$ is cubic and $2$-connected. For $k=0$, the equality above holds since $K_{3,3}$ is a frame of $P_{10}$ and $l(K_{3,3})=0$. 

\noindent
Set $Q:=P^k_{10}$ with $k \geq 1$. 
Let $M$ be an $f$-matching of $Q$. Denote the $2$-connected cubic graph which is homeomorphic to $Q-M$ by $\overline Q(k)$. Suppose that $M$ is chosen in such a way that $l(\overline Q(k))$ is minimal.

\noindent
A subgraph of $\overline Q(k)$ is denoted by $\overline H$, say, and the corresponding subgraph in $Q-M$ and $Q$ by $H$.

\noindent
Let $\overline C$ be a circuit of $\overline Q(k)$ such that $\max\limits_{v \in \overline Q(k)} d_{\overline Q(k)}(\overline C,v)= l(\overline Q(k))$. 
$Q$ contains ten disjoint induced subgraphs isomorphic to $P^{k-1}$. If we contract each of them to a distinct vertex, we obtain $P_{10}$. Hence $C$ does not pass through each of them since otherwise it would imply that $P_{10}$ is hamiltonian. Let us denote one copy of $P^{k-1}$ in $Q$ which is vertex-disjoint with $C$, by $X$.

\noindent
Let $\{v_1,v_2\} \subseteq V(X)$, then Def. \ref{l2} implies, if $v_1$ and $v_2$ are contained in the same copy of $P$, that $p[v_1,W_{k-1}(X)]=p[v_2,W_{k-1}(X)]$.
Therefore and by Lemma \ref{connp2} (2) there is a vertex $x \in V(X)$ which is not matched by $M$ and which satisfies, $p[x, W_{k-1}(X)]= p_{k-1}$, see Def. \ref{l2}. Denote also by $x$ the corresponding vertex in $\overline Q(k)$.

\noindent
Let $\overline \alpha_{x} \subseteq \overline Q(k)$ be a path of length $d(x,\overline C)$ which connects $x$ with $\overline C$.

\noindent
By the definition of $x$, $p(\alpha_x) \geq p_{k-1}$. Since $V(C) \cap V(X) = \emptyset $, $\alpha_x$ passes through at least $p_{k-1}-1$ distinct copies of $P$. For every such copy of $P$, $\overline \alpha_{x}$ contains by Lemma \ref{connp6} at least one vertex. Since $\overline \alpha_{x}$ starts and ends in a vertex of degree $3$ which is not contained in any of these copies of $P$, $|V(\overline \alpha_{x})| \geq p_{k-1}+1$. Thus and by definition of $\overline C$ and 
$\overline \alpha_{x}$, $$l(\overline Q(k)) \geq d(x,\overline C) \geq p_{k-1}\,\,\,\,\,\,\,\,\,\,\,\,\,\,\,\,\,\,\,(4)$$ 

\noindent
Consider $k=1$. By inequality $(4)$, $l(\overline Q(1)) \geq p_0$. Since $p_0=1$ (Lemma \ref{dkak}) and since $P_{10}$ is a frame of $Q$ (Lemma \ref{pkpkpk}) with 
$l(P_{10})=1$ (Corollary \ref{Pkdk}), $l(\overline Q(1))=1$. \\
Consider $k > 1$. By inequality $(4)$ and by Lemma \ref{dkak}, $l(\overline Q(k)) \geq 2^{2k-3}$. Since by Lemma \ref{pkpkpk}, $P^{k-1}_{10}$ is a frame of $Q$ and since by Corollary \ref{Pkdk} $l(P^{k-1}_{10})=2^{2k-3}$, $l(\overline Q(k))= 2^{2k-3}$ which finishes the proof.

\begin{corollary}\label{}
Every $P^k_{10}$, $k \geq 1$ is a counterexample to Conjecture \ref{3cf}. 
\end{corollary}

\begin{corollary}\label{}
For every set $\mathcal S_0$ of $2$-connected graphs with $l_m(\mathcal S_0) \not= \infty$, there is an infinite set $\mathcal S$ of $3$-connected cubic graphs with the following property: for every $G \in \mathcal S$, $\mathcal S_0$ is not an even subdivision-factor of $G$.
\end{corollary}

\noindent
Proof: Replace every element in $\mathcal S_0$ which contains a $2$-valent and a $3$-valent vertex by its homeomorphic cubic graph. Denote this set by $\mathcal T_0$. We observe that if $\mathcal S_0$ is an even subdivision-factor of a cubic graph $H$, say, then $\mathcal T_0$ is also an even subdivision-factor of $H$.
Moreover, $l_m(\mathcal T_0) \leq l_m(\mathcal S_0)$. Set $\mathcal S:=\{P^{k}_{10}\,\,|\,\,2^{2k-3} \,> \,l_m(\mathcal T_0), \,k \geq 2 \}$. 
Theorem \ref{pktk} implies that for every $G \in \mathcal S$, $\mathcal T_0$ is not an even subdivision-factor of $G$. By the above observation, the same holds for $\mathcal S_0$ which finishes the proof.





\footnotesize
\bibliographystyle{plain}

\end{document}